\documentclass[11pt,a4paper]{article}
\usepackage[a4paper]{geometry}
\usepackage{amssymb,latexsym,amsmath,amsfonts,amsthm}
\usepackage{graphicx}
\usepackage{algorithm,algorithmic}
\usepackage{dsfont}
\usepackage{epstopdf}
\usepackage{subfigure}
\usepackage{mathrsfs,enumerate}
\usepackage{tikz}
\usepackage{epsfig}
\usepackage{booktabs}
\usepackage{mathtools}
\usepackage{comment,verbatim}
\usepackage{overpic}
\usepackage{hyperref}
\usepackage{bbm}
\usepackage[toc,page]{appendix}

\newtheorem{theorem}{Theorem}[section]
\newtheorem{lemma}[theorem]{Lemma}

\theoremstyle{definition}

\theoremstyle{definition}
\newtheorem{example}[theorem]{Example}

\newtheorem{remark}[theorem]{Remark}

\numberwithin{equation}{section}

\usepackage{color}

\usepackage[normalem]{ulem}

\begin{document}
\title{Superconvergence points of Hermite spectral interpolation}
\author{Haiyong Wang\footnotemark[1]~\footnotemark[2] ~ and~ Zhimin Zhang\footnotemark[3]}
\date{}
\maketitle
\renewcommand{\thefootnote}{\fnsymbol{footnote}}

\footnotetext[1]{School of Mathematics and Statistics, Huazhong
University of Science and Technology, Wuhan 430074, P. R. China.
E-mail: \texttt{haiyongwang@hust.edu.cn}}

\footnotetext[2]{Hubei Key Laboratory of Engineering Modeling and
Scientific Computing, Huazhong University of Science and Technology,
Wuhan 430074, P. R. China}

\footnotetext[3]{Department of Mathematics, Wayne State University, Detroit, MI 48202, USA}

\begin{abstract}
Hermite spectral method plays an important role in the numerical simulation of various partial differential equations (PDEs) on unbounded domains. In this work, we study the superconvergence properties of Hermite spectral interpolation, i.e., interpolation at the zeros of Hermite polynomials in the space spanned by Hermite functions. We identify the points at which the convergence rates of the first- and second-order derivatives of the interpolant converge faster. We further extend the analysis to the Hermite spectral collocation method in solving differential equations and identify the superconvergence points both for function and derivative values. Numerical examples are provided to confirm the analysis of superconvergence points.
\end{abstract}


\vspace{0.05in}


\section{Introduction}\label{sec:introduction}
Spectral methods are built upon spectral approximations using classical orthogonal polynomials and they are widely used in the numerical simulation of various differential or integral equations (see \cite{Boyd2000,Canuto2006,Shen2011}). It is well known that spectral methods have the remarkable advantage that they exhibit so-called spectral convergence, i.e., their convergence depends solely on the regularity of the underlying functions. In particular, the convergence can be improved to exponential if the underlying functions are analytic (see, e.g., \cite{Bern1912,Reddy2005,Tadmor1986,Wang2012,Wang2021,Wang2024,Wang2025a,Wang2014,Xiang2010,Xie2013,Zhao2013}).

The superconvergence phenomenon of $h$-version methods has received considerable attention in the past few decades, including the finite element method for PDEs \cite{Wahlbin1995}, the collocation method for Volterra integral equations \cite{Brunner2004}, and the discontinuous Galerkin method for hyperbolic equations \cite{Cao2014,Cao2015}. As for spectral method, however, only a few studies can be found in the literature (see, e.g., \cite{Wang2025b,Wang2014,Zhang2005,Zhang2012,Zhao2016}). Let $\Omega=[-1,1]$ and let $p_n$ denote some spectral approximation of degree $n$ to $f$. For each $k\in\mathbb{N}$, the superconvergence analysis of $p_n$ is to find a set of points $\{y_j\}\subset\Omega$ that satisfy
\[
n^{\alpha} |(f - p_n)^{(k)}(y_j)| \leq C \|(f - p_n)^{(k)}\|_{L^{\infty}(\Omega)},
\]
for some $C>0$ and $\alpha>0$. If $f$ is analytic in a region containing $\Omega$, the superconvergence points of the Chebyshev spectral interpolation were first studied in \cite{Zhang2012}. It was shown that the convergence rate of the first- and second-order derivatives of the Chebyshev spectral interpolation can be improved by a factor of $n$ or $n^2$ at the superconvergence points. The analysis was subsequently extended to the case of Jacobi spectral interpolation in \cite{Wang2014} and to spectral interpolation involving fractional derivatives in \cite{Zhao2016}, with similar superconvergence results derived. If $f$ is only a   differentiable function on $\Omega$, however, the superconvergence analysis is much more involved. More specifically, it was shown in \cite{Wang2014} that the accuracy of Jacobi spectral interpolation can still be improved at the superconvergence points of analytic cases, but there is no gain in the order of convergence. More recently, for functions with an algebraic singularity, the superconvergence points of the Jacobi projection were studied in \cite{Wang2025b}. It was shown that the derivatives of the Jacobi projection superconverge at those superconvergence points when they are bounded away from the singularity.

In this paper, we study the superconvergence points of Hermite spectral interpolation, which uses Hermite functions and interpolates at the zeros of Hermite polynomials. Note that most existing superconvergence analyses focus on spectral methods on finite intervals, and very little is known about the superconvergence phenomenon of spectral methods in unbounded intervals.
Building upon the recent studies on the convergence analysis of Hermite spectral approximation in \cite{Wang2025a}, we identify the superconvergence points of Hermite spectral interpolation, at which the derivatives can be improved by at least a factor of $n^{1/2}$, where $n$ is the degree of Hermite spectral interpolation. We extend the superconvergence analysis to the Hermite spectral collocation method for two ODE models and identify the superconvergence points at which the Hermite spectral collocation method converges faster for both function and derivative values. Further, we discuss the post-processing of the Hermite spectral collocation method based on the superconvergence analysis. We show that the accuracy of Hermite spectral collocation methods can be further improved using their information at the superconvergence points.

The rest of this paper is organized as follows. In the next section, we recall some properties of Hermite polynomials and functions. In section \ref{sec:Super} we analyze the error of Hermite spectral interpolation and identify the superconvergence points for the first- and second-order derivatives. In section \ref{sec:HermSpecColl} we extend the analysis to Hermite spectral collocation method for two
second order linear ODEs with constant and variable coefficients and and identify the superconvergence points at which Hermite spectral collocation method has the superconvergence properties for function as well as derivative values. In section \ref{sec:Post} we extend the analysis to the post-processing of Hermite spectral collocation method. We conclude our results in section \ref{sec:Conclusion}.

\section{Hermite polynomials and functions}\label{sec:HermPoly}
In this section, we review some basic properties of Hermite
polynomials and functions. We refer to \cite{Szego1939} for a more thorough treatment. Throughout the paper, we denote by $\mathbb{N}_0$ the set of nonnegative integers.

For $n\in\mathbb{N}_0$, the Hermite polynomials are defined by
\begin{equation}\label{def:HermPoly}
H_n(x) = (-1)^n \mathrm{e}^{x^2} \frac{\mathrm{d}^n}{\mathrm{d}x^n}\mathrm{e}^{-x^2},
\end{equation}
and it is well-known that
\begin{equation}\label{eq:ortho}
\int_{-\infty}^{\infty} H_n(x) H_m(x) \mathrm{e}^{-x^2} \mathrm{d}x = \gamma_{n} \delta_{n,m},
\end{equation}
where $\delta_{n,m}$ is the Kronecker delta and $\gamma_{n}=2^n n!\sqrt{\pi} $. Below we list some properties of Hermite polynomials:
\begin{itemize}
  \item Recurrence relation:
  \begin{equation}\label{eq:recur}
  H_{n+1}(x) = 2xH_n(x) - 2nH_{n-1}(x), \quad n\geq1,
  \end{equation}
  where $H_0(x) = 1$ and $H_1(x) = 2x$.

  \item Symmetry:
  \begin{equation}\label{eq:HermSymm}
  H_n(-x) = (-1)^n H_n(x), \quad n\geq0,
  \end{equation}
  and thus $H_n(x)$ is an even function when $n$ is even and $H_n(x)$ is an odd function when $n$ is odd.

  \item Derivative:
  \begin{equation}\label{eq:HermD}
  H'_n(x)=2nH_{n-1}(x), \quad n\geq1.
  \end{equation}
\end{itemize}
For $n\in\mathbb{N}_0$, the Hermite functions are defined by
\begin{equation}\label{def:HermFunc}
\psi_n(x) = \mathrm{e}^{-x^2/2}\frac{H_n(x)}{\sqrt{\gamma_n}} .
\end{equation}
It is well-known that $\{\psi_n\}_{n=0}^{\infty}$ form a complete orthonormal basis for the Hilbert space $L^2(\mathbb{R})$. Below we list some important properties of Hermite functions:
\begin{itemize}
\item They satisfy the following inequality
\begin{equation}\label{eq:PsiBound}
|\psi_n(x)| \leq \frac{1}{\pi^{1/4}}, \quad  x\in\mathbb{R}.
\end{equation}
Moreover, this inequality is sharp in the sense that the upper bound can be attained when $n=0$ and $x=0$ \cite{Indritz1961}.

\item By \eqref{eq:recur}, the recurrence relation of $\{\psi_n\}$ is
\begin{equation}\label{eq:recur2}
x\psi_{n}(x) = \sqrt{\frac{n}{2}} \psi_{n-1}(x) + \sqrt{\frac{n+1}{2}} \psi_{n+1}(x), % \quad n\geq0,
\end{equation}
where $\psi_0(x)=\mathrm{e}^{-x^2/2}/\pi^{1/4}$ and $\psi_1(x) = x\mathrm{e}^{-x^2/2}\sqrt{2}/\pi^{1/4}$.

\item The derivatives of Hermite functions satisfy
\begin{align}\label{eq:Deriv}
\psi'_n(x) &= \sqrt{\frac{n}{2}} \psi_{n-1}(x) - \sqrt{\frac{n+1}{2}} \psi_{n+1}(x).  % \quad n\geq0.
\end{align}

\item Hermite functions are the eigenfunctions of the harmonic oscillator:
\begin{equation}\label{eq:harm}
\left(- \frac{\mathrm{d}^2}{\mathrm{d}x^2} + x^2 \right) \psi_n(x) = (2n+1) \psi_n(x).  %, \quad n\in\mathbb{N}_0.
\end{equation}
\end{itemize}

\section{Superconvergence points of Hermite spectral interpolations}\label{sec:Super}
In this section we consider Hermtie spectral interpolation method. Let $\{x_j\}_{j=0}^{n}$ be the zeros of $H_{n+1}(x)$ and we assume that they are arranged in ascending order, i.e., $-\infty<x_0<x_1<\cdots<x_n<\infty$. By the symmetry relation \eqref{eq:HermSymm}, we know that $x_j=-x_{n-j}$ for $j=0,\ldots,n$. Let $\mathbb{H}_{n}$ denote the space spanned by $\{\psi_k\}_{k=0}^{n}$, i.e., $\mathbb{H}_n=\mathrm{span}\{\psi_k\}_{k=0}^{n}$, and let $h_{n}\in\mathbb{H}_{n}$ be the unique function which interpolates $f(x)$ at the points
$\{x_j\}_{j=0}^{n}$, i.e.,
\begin{equation}\label{eq:hn}
h_{n}(x_j) = f(x_j), \qquad j=0,\ldots,n.
\end{equation}
Inspired by \cite{Zhang2012}, we consider the superconvergence points of the interpolant $h_{n}$ and find the set of points $\{y_j\}\subset\mathbb{R}$ such that
\begin{equation}
n^{\alpha} |(f - h_n)^{(k)}(y_j) | \leq  C \|(f - h_n)^{(k)}\|_{L^{\infty}(\mathbb{R})},
\end{equation}
where $\alpha>0$ and $C$ is some positive constant, for $k=1,2$.

We start with a contour integral representation of the remainder of $h_n(x)$, which was recently proved in \cite[Theorem 4.5]{Wang2025a}. Before proceeding, we introduce the infinite strip in the complex plane
\begin{equation}\label{def:Strip}
\mathcal{S}_{\rho} := \big\{z\in\mathbb{C}: ~ \Im(z) \in [-\rho,\rho] \big\},
\end{equation}
where $\rho\in(0,\infty)$, and denote by $\partial{\mathcal{S}}_{\rho}$ the boundary of $\mathcal{S}_{\rho}$. It is known that the convergence domain of Hermite approximation in the complex plane can be characterized by an infinite strip for analytic functions (see, e.g., \cite[Section~9.2]{Szego1939}). In what follows, the orientation of contour integrals along $\partial{\mathcal{S}}_{\rho}$ is always taken from left to right when $\Im(z)=-\rho$ and from right to left when $\Im(z)=\rho$. Throughout the paper, we denote by $\mathcal{K}$ a generic positive constant and by $\mathrm{i}$ the imaginary unit.
\begin{lemma}\label{lem:HermCont}
If $f$ is analytic in the strip $\mathcal{S}_{\rho}$ for
some $\rho>0$ and $|\mathrm{e}^{z^2/2}f(z)|\leq\mathcal{K}|z|^{\sigma}$ for some
$\sigma\in\mathbb{R}$ as $|z|\rightarrow\infty$ within the strip and if
\begin{equation}
\widehat{V} := \int_{\partial\mathcal{S}_{\rho}} |\mathrm{e}^{z^2/2} f(z)| |\mathrm{d}z| < \infty. \nonumber
\end{equation}
Then, for $n\geq\max\{\lfloor\sigma\rfloor,0\}$,
\begin{equation}
f(x) - h_{n}(x) =
\frac{1}{2\pi\mathrm{i}} \int_{\partial{\mathcal{S}}_{\rho}}
\frac{\psi_{n+1}(x)f(z)}{\psi_{n+1}(z) (z-x)} \mathrm{d}z. \nonumber
\end{equation}
\end{lemma}

Now we consider the maximum error estimates of the first- and second-order spectral differentiations using $h_n$. Differentiating the contour integral representation of $f(x) - h_{n}(x)$ in Lemma \ref{lem:HermCont} once and twice, we obtain
\begin{align}\label{eq:FirD}
(f - h_{n}){'}(x) =
\frac{1}{2\pi\mathrm{i}} \int_{\partial{\mathcal{S}}_{\rho}}
\frac{f(z)}{\psi_{n+1}(z)} \left(\frac{\psi'_{n+1}(x)}{z-x} + \frac{\psi_{n+1}(x)}{(z-x)^2} \right) \mathrm{d}z,
\end{align}
and
\begin{align}\label{eq:SecD}
(f - h_{n}){''}(x) =
\frac{1}{2\pi\mathrm{i}} \int_{\partial{\mathcal{S}}_{\rho}}
\frac{f(z)}{\psi_{n+1}(z)} \left(\frac{\psi''_{n+1}(x)}{z-x} + \frac{2\psi'_{n+1}(x)}{(z-x)^2} + \frac{2\psi_{n+1}(x)}{(z-x)^3} \right) \mathrm{d}z.
\end{align}
Note that the interchanges of the order of differentiation and integration in \eqref{eq:FirD} and \eqref{eq:SecD} can be easily justified by the assumptions of $f(z)$ and the dominated convergence theorem. Clearly, to establish bounds for $(f - h_{n}){'}(x)$ and $(f - h_{n}){''}(x)$, it is necessary to estimate the maximum norm of the derivatives of $\psi_{n}(x)$ for $x\in\mathbb{R}$. We establish the result in the following lemma.
\begin{lemma}\label{lem:PsiD}
Let $k\in\mathbb{N}_0$ and $n\in\mathbb{N}_0$. For $n\geq1$, it holds that
\begin{equation}\label{eq:PsiBound1}
\|\psi_n^{(k)}\|_{L^{\infty}(\mathbb{R})} \leq \mathcal{C}_k \left\{
           \begin{array}{ll}
            n^{-1/12},                   & \hbox{$k=0$,} \\[2ex]
           {\displaystyle  n^{k/2-1/4}}, & \hbox{$k\geq1$,}
           \end{array}  \right.
\end{equation}
where $\mathcal{C}_k$ are some positive constant independent of $n$. For $n=0$, it holds that
\begin{equation}\label{eq:PsiBound2}
\|\psi_0^{(k)}\|_{L^{\infty}(\mathbb{R})} \leq \left\{
           \begin{array}{ll}
            \pi^{-1/4},                   & \hbox{$k=0$,} \\[2ex]
           {\displaystyle 2^{-1/2}\mathcal{C}_{k-1} }, & \hbox{$k\geq1$.}
           \end{array}  \right.
\end{equation}
\end{lemma}
\begin{proof}
We first prove \eqref{eq:PsiBound1}. The case $k=0$ follows from \cite[Equation (30)]{Hille1926}. For $k=1$, we obtain from \eqref{eq:Deriv} that
\begin{align}
\psi'_n(x) &= \sqrt{\frac{n}{2}} \psi_{n-1}(x) - \sqrt{\frac{n+1}{2}} \psi_{n+1}(x) \nonumber \\
&= -\sqrt{\frac{n+1}{2}} (\psi_{n+1}(x) - \psi_{n-1}(x)) - \frac{\sqrt{n+1}-\sqrt{n}}{\sqrt{2}} \psi_{n-1}(x). \nonumber
\end{align}
As $n\rightarrow\infty$, it is easily seen that the second term on the right-hand side behaves like $O(n^{-7/12})$. Moreover, by \cite[Lemma 1]{Markett1984} we know that $\|\psi_{n+1}-\psi_{n-1}\|_{L^{\infty}(\mathbb{R})}=O(n^{-1/4})$, and thus the first term on the right-hand side behaves like $O(n^{1/4})$. Combining these two estimates we conclude that $\|\psi'_n\|_{L^{\infty}(\mathbb{R})}=O(n^{1/4})$. This proves $k=1$. As for $k\geq2$, by \eqref{eq:Deriv} again we have
\begin{align}
\psi_n^{(k)}(x) &= \sqrt{\frac{n}{2}} \psi_{n-1}^{(k-1)}(x) - \sqrt{\frac{n+1}{2}} \psi_{n+1}^{(k-1)}(x). \nonumber
\end{align}
The desired result \eqref{eq:PsiBound1} then follows by induction on $k$. We next consider \eqref{eq:PsiBound2}. The case $k=0$ follows immediately from \eqref{eq:PsiBound}. For $k\geq1$, note that $\psi_0^{(k)}(x)=-\psi_{1}^{(k-1)}(x)/\sqrt{2}$, the desired result \eqref{eq:PsiBound2} then follows by combining this with \eqref{eq:PsiBound1}. This ends the proof.
\end{proof}

\begin{figure}[htbp]
\centering
\includegraphics[width=0.7\textwidth,height=0.55\textwidth]{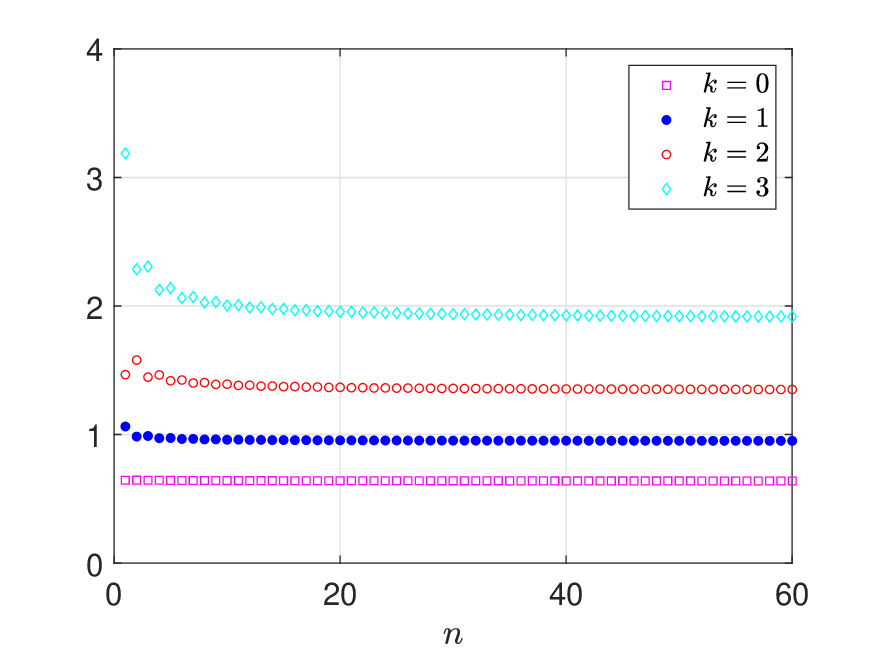}
\caption{Plot of $\|\psi_n\|_{L^{\infty}(\mathbb{R})} n^{1/12}$ and $\|\psi_n^{(k)}\|_{L^{\infty}(\mathbb{R})} n^{1/4-k/2}$ for $k=1,2,3$ as functions of $n$. Here $n$ ranges from $1$ to $60$.}\label{fig:HermMax}
\end{figure}

In Figure \ref{fig:HermMax} we plot $\|\psi_n\|_{L^{\infty}(\mathbb{R})} n^{1/12}$ and $\|\psi_n^{(k)}\|_{L^{\infty}(\mathbb{R})} n^{1/4-k/2}$ for $k=1,2,3$ as functions of $n$. We observe that they either tend to some constants or decay rather slowly as $n$ increases, and thus we can conclude that \eqref{eq:PsiBound1} is sharp or nearly sharp.

\begin{remark}
Numerical calculations show that $\|\psi_n\|_{L^{\infty}(\mathbb{R})} n^{1/12}$ and $\|\psi''_n\|_{L^{\infty}(\mathbb{R})} n^{-3/4}$ attain their maximum when $n=2$ and $\|\psi'_n\|_{L^{\infty}(\mathbb{R})} n^{-1/4}$ and $\|\psi'''_n\|_{L^{\infty}(\mathbb{R})} n^{-5/4}$ attain their maximum when $n=1$. Thus, the optimal constants of $\{\mathcal{C}_k\}_{k=0}^{3}$, i.e., the smallest constants such that \eqref{eq:PsiBound1} holds for all $n\in\mathbb{N}$, can be calculated as
\begin{align}
\mathcal{C}_0 &= \|\psi_2\|_{L^{\infty}(\mathbb{R})} 2^{1/12} = \frac{2^{19/12} \mathrm{e}^{-5/4}}{\pi^{1/4}} \approx 0.644874576859960, \nonumber \\
\mathcal{C}_1 &= \|\psi'_1\|_{L^{\infty}(\mathbb{R})} = \frac{\sqrt{2}}{\pi^{1/4}} \approx  1.062251932027197,   \nonumber \\
\mathcal{C}_2 &= \|\psi''_2\|_{L^{\infty}(\mathbb{R})} 2^{-3/4} = \frac{5}{2^{5/4}\pi^{1/4}} \approx 1.579046944365162, \nonumber \\
\mathcal{C}_3 &= \|\psi'''_1\|_{L^{\infty}(\mathbb{R})} = \frac{3\sqrt{2}}{\pi^{1/4}} \approx 3.186755796081591. \nonumber
\end{align}
\end{remark}

\begin{remark}
The following inequalities were proved in \cite[Lemma~8.3]{Comte2018} when dealing with inverse problems of nonparametric statistics,
\begin{equation}\label{eq:Comte}
\|\psi'_n\|_{L^{\infty}(\mathbb{R})} \leq \mathcal{K}_1 (n+1)^{5/12}, \quad
\|(x\psi_n)'\|_{L^{\infty}(\mathbb{R})} \leq \mathcal{K}_2 (n+1)^{11/12},
\end{equation}
where $\mathcal{K}_1$ and $\mathcal{K}_2$ are some positive constants. For large $n$, our result in \eqref{eq:PsiBound1} gives $\|\psi'_n\|_{L^{\infty}(\mathbb{R})}=O(n^{1/4})$, which is better than the former inequality in \eqref{eq:Comte}. Moreover, by combining \eqref{eq:recur2} and \eqref{eq:PsiBound1}, we obtain for $n\geq2$ that
\begin{align}
\|(x\psi_n)'\|_{L^{\infty}(\mathbb{R})} &\leq \sqrt{\frac{n+1}{2}} \|\psi'_{n+1}\|_{L^{\infty}(\mathbb{R})} + \sqrt{\frac{n}{2}} \|\psi'_{n-1}\|_{L^{\infty}(\mathbb{R})} \leq \mathcal{C}_1 \sqrt{2} (n+1)^{3/4},
\end{align}
which is also better than the latter inequality in \eqref{eq:Comte}.
\end{remark}

\begin{remark}\label{rk:EstH}
For large $n$, the estimate of $|\psi_n(x)|$ will be different when $x\in[-\varrho\xi,\varrho\xi]$, where $\xi=(2n+1)^{1/2}$ is the turning point of $\psi_n(x)$ and $\varrho\in(0,1)$. In this case, from \cite[Equation~(29)]{Hille1926} we know that $|\psi_n(x)|=O(n^{-1/4})$ (see also \cite[Equation~(8.22.12)]{Szego1939}). Indeed, the location where $|\psi_n(x)|$ attains its maximum for $x\in\mathbb{R}$ will approach to the turning point $\xi$ as $n$ increases (see \cite[Equation (8.22.14)]{Szego1939}).
\end{remark}

With the above two lemmas, we now establish maximum error estimate for the Hermite spectral differentiation.
\begin{theorem}
If $f$ is analytic in the strip $\mathcal{S}_{\rho}$ for
some $\rho>0$ and $|\mathrm{e}^{z^2/2}f(z)|\leq\mathcal{K}|z|^{\sigma}$ for some
$\sigma\in\mathbb{R}$ as $|z|\rightarrow\infty$ within the strip and if
$\widehat{V} < \infty$. Then, for $n\geq\max\{\lfloor\sigma\rfloor,1\}$ and $m\in\mathbb{N}_0$,
\begin{align}\label{eq:MaxErr}
\|(f - h_{n})^{(m)}\|_{L^{\infty}(\mathbb{R})} \leq \mathcal{K}_m \left\{
           \begin{array}{ll}
           {\displaystyle  n^{1/6} \mathrm{e}^{-\rho\sqrt{2n}} },  & \hbox{$m=0$,} \\[2ex]
           {\displaystyle  n^{m/2} \mathrm{e}^{-\rho\sqrt{2n}} },  & \hbox{$m\geq1$,}
           \end{array}  \right.
\end{align}
where $\mathcal{K}_m=\mathrm{e}^{\rho^2/2}\widehat{V}\mathcal{C}_m/(\pi^{1/2}2^{1/4}\rho)(1+o(1))$ as $n\rightarrow\infty$ and $\mathcal{C}_m$'s are the positive constants defined in \eqref{eq:PsiBound1}.
\end{theorem}
\begin{proof}
We sketch the proof of the case $m=1$ and the proof of the other cases is similar.
Note that $\min|z-x|=\rho$ for $x\in\mathbb{R}$ and $z\in\partial\mathcal{S}_{\rho}$, by \eqref{eq:FirD} and Lemma \ref{lem:PsiD} we have
\begin{align}
|(f - h_{n}){'}(x)| &\leq
\frac{1}{2\pi} \int_{\partial{\mathcal{S}}_{\rho}}
\frac{|f(z)|}{|\psi_{n+1}(z)|} \left( \frac{|\psi'_{n+1}(x)|}{|z-x|} + \frac{|\psi_{n+1}(x)|}{|z-x|^2} \right) |\mathrm{d}z|  \nonumber \\
&\leq \frac{1}{2\pi} \left(\frac{\|\psi'_{n+1}\|_{L^{\infty}(\mathbb{R})}}{\rho} + \frac{\|\psi_{n+1}\|_{L^{\infty}(\mathbb{R})}}{\rho^2} \right) \int_{\partial{\mathcal{S}}_{\rho}}
\frac{|f(z)|}{|\psi_{n+1}(z)|} |\mathrm{d}z|  \nonumber \\
&\leq \frac{1}{2\pi} \left(\frac{\|\psi'_{n+1}\|_{L^{\infty}(\mathbb{R})}}{\rho} + \frac{\|\psi_{n+1}\|_{L^{\infty}(\mathbb{R})}}{\rho^2} \right) \frac{\widehat{V} \sqrt{\gamma_{n+1}} }{\displaystyle \min_{z\in\partial{\mathcal{S}}_{\rho}}|H_{n+1}(z)|}  \nonumber \\
&\leq \frac{1}{2\pi} \left( \frac{\mathcal{C}_1(n+1)^{1/4}}{\rho} + \frac{\mathcal{C}_0(n+1)^{-1/12}}{\rho^2} \right) \frac{\widehat{V} \sqrt{\gamma_{n+1}} }{\displaystyle \min_{z\in\partial{\mathcal{S}}_{\rho}}|H_{n+1}(z)|}.  \nonumber
\end{align}
Furthermore, by \cite[Lemma 4.2]{Wang2025a} we know for $k\in\mathbb{N}$ that
\begin{align}
\min_{z\in\partial\mathcal{S}_{\rho}}|H_{k}(z)|  &= |H_{k}(\mathrm{i}\rho)| \left\{
           \begin{array}{ll}
             1,              & \hbox{$k$ odd,} \\[2ex]
           {\displaystyle  1 + O\left(\mathrm{e}^{-2\rho\sqrt{2k}}\right) }, & \hbox{$k$ even,}
           \end{array}  \right. \nonumber \\[2ex]
           &= \frac{\mathrm{e}^{-\rho^2/2}}{2} \frac{\Gamma(k+1)}{\Gamma(k/2+1)} \mathrm{e}^{\rho\sqrt{2k}} \left(1 + O\left(k^{-1/2}\right)\right),  \nonumber
\end{align}
as $k\rightarrow\infty$, and we have used \cite[Theorem~8.22.7]{Szego1939} in the last step. Combining the above two results yields
\begin{align}
|(f - h_{n}){'}(x)| &\leq \mathcal{K} \frac{\Gamma((n+3)/2) 2^{(n+1)/2}}{\sqrt{\Gamma(n+2)}} n^{1/4}  \mathrm{e}^{-\rho\sqrt{2n}} , \nonumber
\end{align}
where $\mathcal{K}=\mathrm{e}^{\rho^2/2}\widehat{V} \mathcal{C}_1/(\pi^{3/4}\rho)(1 + O(n^{-1/3}))$ as $n\rightarrow\infty$. By the duplication formula \cite[Equation~(5.5.5)]{Olver2010}, we have
\[
\Gamma(n+2) = \frac{2^{n+1}}{\sqrt{\pi}} \Gamma\left(\frac{n+2}{2}\right) \Gamma\left(\frac{n+3}{2}\right),
\]
and thus
\begin{align}
\frac{\Gamma((n+3)/2) 2^{(n+1)/2}}{\sqrt{\Gamma(n+2)}} &= \pi^{1/4} \sqrt{\frac{\Gamma((n+3)/2)}{\Gamma((n+2)/2)}} = \left(\frac{n\pi}{2}\right)^{1/4} \left(1 + O\left(n^{-1}\right)\right),   \nonumber
\end{align}
as $n\rightarrow\infty$, where we have used the asymptotic of the ratio of gamma functions \cite[Equation~(5.11.13)]{Olver2010}. The desired result follows immediately and this ends the proof.
\end{proof}

Let $\{\tau_j\}_{j=0}^{n+1}$ denote the zeros of $\psi'_{n+1}(x)$. Note that the first term inside the parentheses of the equation \eqref{eq:FirD} vanishes at the points $\{\tau_j\}_{j=0}^{n+1}$, by Lemma \ref{lem:PsiD} we have
\begin{equation}
|(f - h_n){'}(\tau_j) | \leq \mathcal{K} n^{1/6} \mathrm{e}^{-\rho\sqrt{2n}},
\end{equation}
where $\mathcal{K}=\mathrm{e}^{\rho^2/2}\widehat{V}\mathcal{C}_0/(\pi^{1/2}2^{1/4}\rho^2)(1+o(1))$ as $n\rightarrow\infty$. We see that the error bound at the points $\{\tau_j\}_{j=0}^{n+1}$ is smaller than the maximum error bound in \eqref{eq:MaxErr} by a factor $n^{1/3}$. We point out that this factor can be further improved at those points $\{\tau_j\}\subseteq[-\varrho\xi,\varrho\xi]$, where $\varrho$ and $\xi$ are defined as in Remark \ref{rk:EstH}. In this case, note that the estimate $|\psi_{n+1}(\tau_j)|\leq \|\psi_{n+1}\|_{L^{\infty}(\mathbb{R})}=O(n^{-1/12})$ can be replaced by the more sharper estimate $|\psi_{n+1}(\tau_j)|=O(n^{-1/4})$, we then obtain a smaller bound
\begin{equation}
|(f - h_n){'}(\tau_j) | \leq \mathcal{K} \mathrm{e}^{-\rho\sqrt{2n}}.
\end{equation}
Thus, the error bound at those points $\{\tau_j\}\subseteq[-\varrho\xi,\varrho\xi]$ is smaller than the maximum error bound in \eqref{eq:MaxErr} by a factor $n^{1/2}$. Indeed, as shown by the examples given below, the maximum error of the Hermite spectral differentiation is often attained near the origin for rapidly decaying functions, and we can therefore expect that the factor gained at the points $\{\tau_j\}_{j=0}^{n+1}$ will be $n^{1/2}$. Furthermore, let $\{\eta_j\}_{j=0}^{n+2}$ denote the zeros of $\psi''_{n+1}(x)$. Note that the first term inside the parentheses of the equation \eqref{eq:SecD} vanishes at the points $\{\eta_j\}_{j=0}^{n+2}$, we have
\begin{equation}
|(f - h_n){''}(\eta_j) | \leq \mathcal{K} n^{1/2} \mathrm{e}^{-\rho\sqrt{2n}},
\end{equation}
where $\mathcal{K}=2^{3/4}\mathrm{e}^{\rho^2/2}\widehat{V}C_1/(\pi^{1/2}\rho^2)(1+o(1))$ as $n\rightarrow\infty$. We see that the error bound at the points $\{\eta_j\}_{j=0}^{n+2}$ is smaller than the maximum error bound in \eqref{eq:MaxErr} by a factor $n^{1/2}$.

\begin{figure}[ht]
\centering
\includegraphics[width=.48\textwidth,height=.4\textwidth]{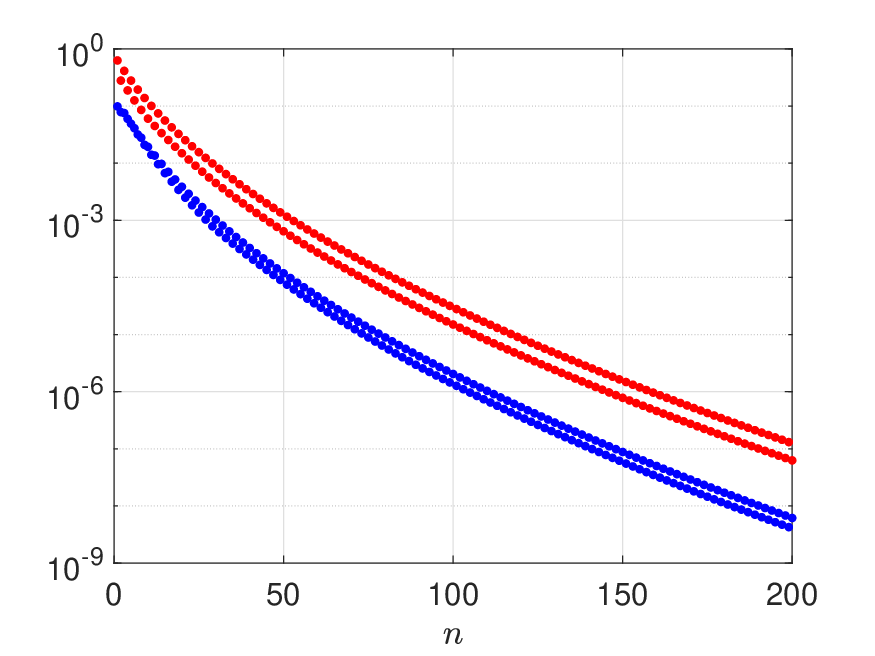}~
\includegraphics[width=.48\textwidth,height=.4\textwidth]{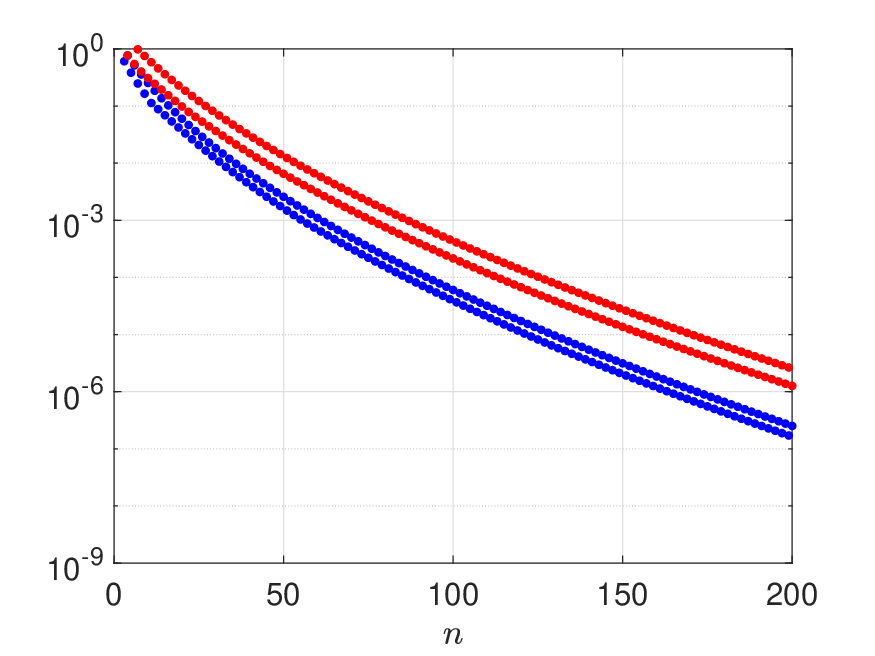}\\
\includegraphics[width=.48\textwidth,height=.4\textwidth]{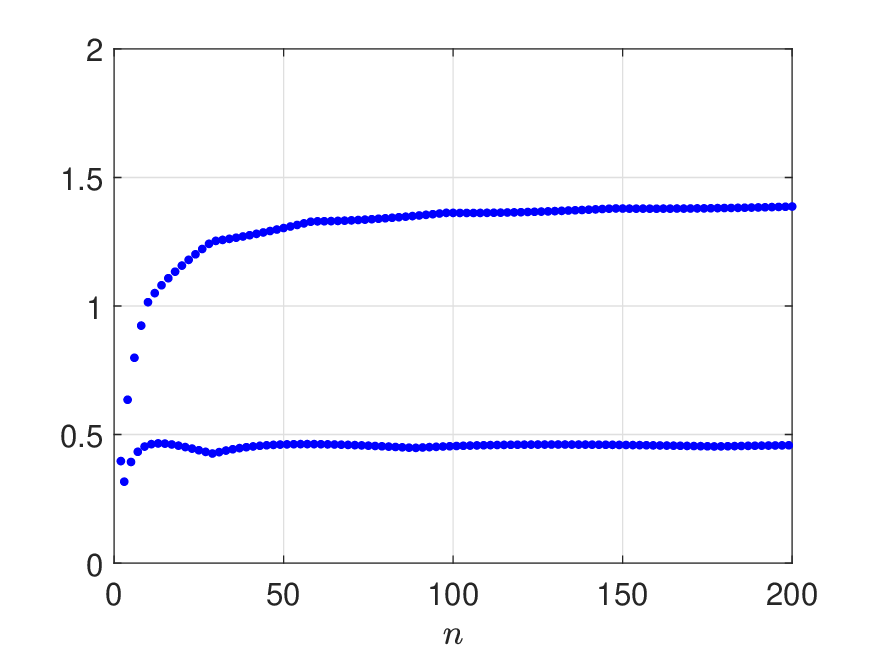}~
\includegraphics[width=.48\textwidth,height=.4\textwidth]{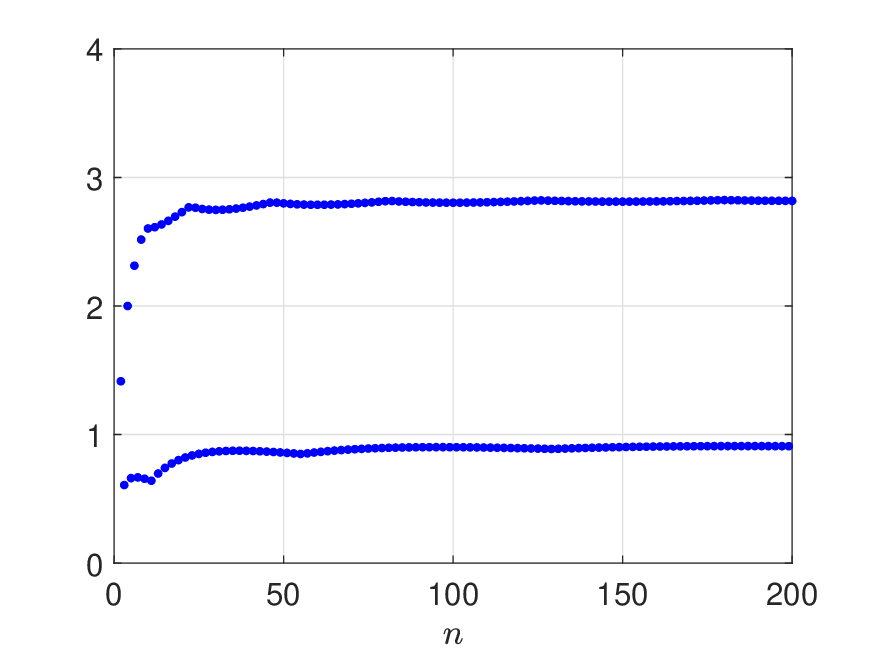}
\caption{Top row: The left panel plots $\max_{j=0,\ldots,n+1}|(f - h_n){'}(\tau_j)|$ (blue) and $\|(f-h_n){'}\|_{\infty}$ (red) and the right panel plots $\max_{j=0,\ldots,n+2}|(f - h_n){''}(\eta_j)|$ (blue) and $\|(f-h_n){''}\|_{\infty}$ (red). Bottom row: Plot of $n^{1/2}\mathcal{R}_1(n)$ (left) and $n^{1/2}\mathcal{R}_2(n)$ (right) as functions of $n$. Here $f(x)=\mathrm{e}^{-x^2/2}/(x^2+1)$ and $n$ ranges from $1$ to $200$.}\label{fig:MaxError}
\end{figure}

Regarding the supervonvergence points of Hermite and Chebyshev spectral interpolations, we make the following remarks.
\begin{remark}
For Chebyshev spectral interpolation, such as polynomial interpolation at the zeros of $T_{n+1}(x)$, it has been shown in \cite{Zhang2012} that the supervonvergence points is the zeros of $T'_{n+1}(x)$ for the first order spectral differentiation and the zeros of $T''_{n+1}(x)$ for the second order spectral differentiation. A factor $n^2$ or $n$ can be gained for the first- and second-order Chebyshev spectral differentiations at the superconvergence points. As for the case of Hermite spectral interpolation, however, the factor gained at the superconvergence points is $n^{1/2}$ (see Figure \ref{fig:MaxError} for an illustration). Note that both methods are used in different contexts (i.e., $[-1,1]$ versus $(-\infty,\infty)$) and their convergence behaviors for analytic functions are also different (i.e., exponential versus root-exponential).
\end{remark}

\begin{remark}
We see from \eqref{eq:harm} that $\psi''_{n+1}(x)=(x^2-2n-3)\psi_{n+1}(x)$, and thus the zeros of $\psi''_{n+1}(x)$ satisfy $\{\eta_j\}_{j=0}^{n+2}=\{x_j\}_{j=0}^{n}\cup\{\pm(2n+3)^{1/2}\}$. This means that, for Hermite spectral interpolation, the second derivative actually superconverges at the interpolation points $\{x_j\}_{j=0}^{n}$. Further, for Hermite spectral collocation method, this justifies the accuracy advantage of the second order Hermite differentiation matrix (see the matrix $D$ in section \ref{sec:HermSpecColl}) in approximating the second derivative values at the interpolation points.
\end{remark}

\begin{figure}[htbp]
\centering
\includegraphics[width=0.7\textwidth,height=0.55\textwidth]{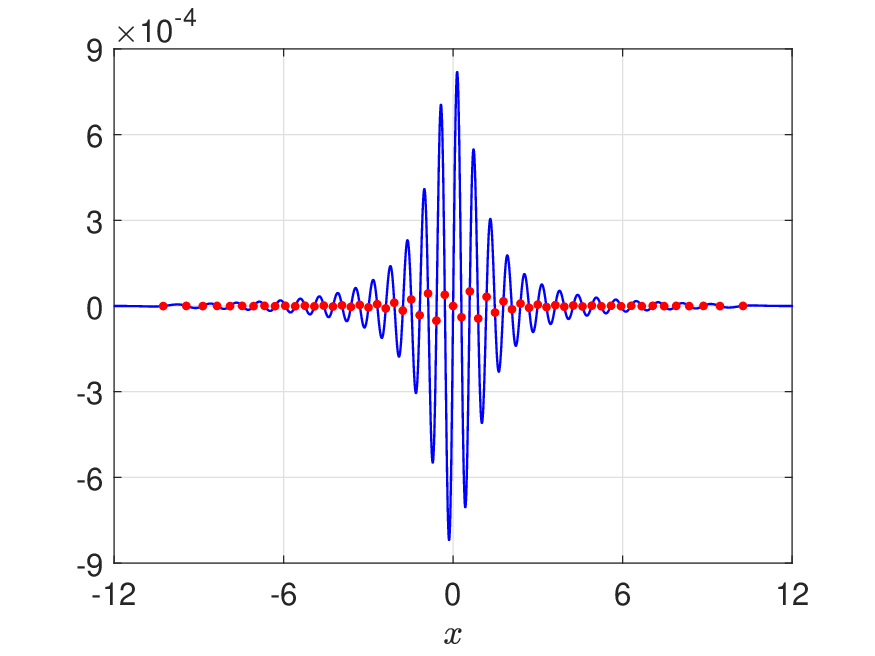}
\caption{Plot of $(f-h_n){'}(x)$ for $f(x)=\mathrm{e}^{-x^2/2}/(x^2+1)$ and $n=55$ and the points are the errors at the superconvergence points $\{\tau_j\}_{j=0}^{n+1}$.}\label{fig:Super1}
\end{figure}

\begin{figure}[htbp]
\centering
\includegraphics[width=0.7\textwidth,height=0.55\textwidth]{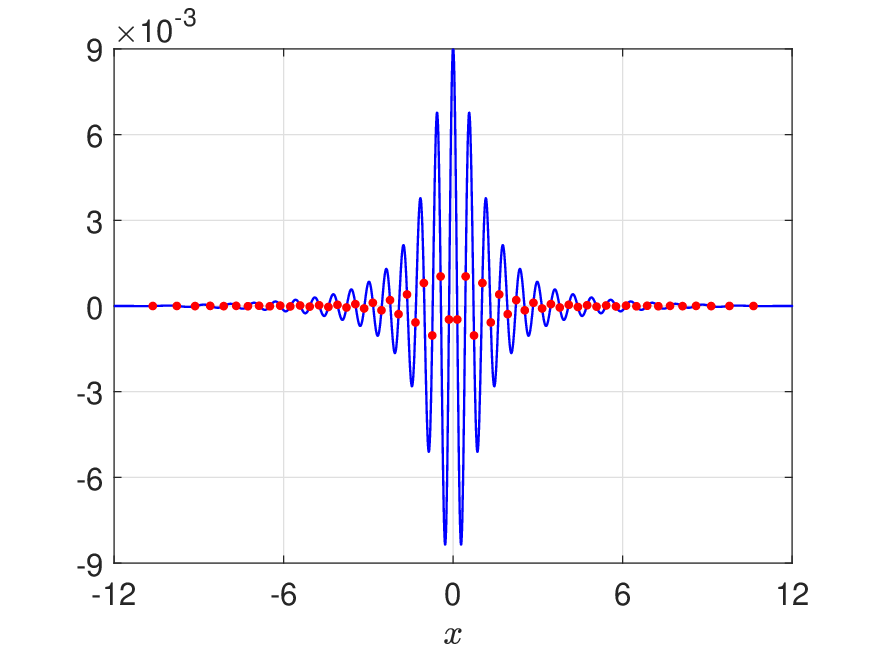}
\caption{Plot of $(f-h_n){''}(x)$ for $f(x)=\mathrm{e}^{-x^2/2}/(x^2+1)$ and $n=55$ and the points are the errors at the superconvergence points $\{\eta_j\}_{j=0}^{n+2}$.}\label{fig:Super2}
\end{figure}

\begin{example}
We consider the function
\[
f(x) = \frac{\mathrm{e}^{-x^2/2}}{x^2+1},
\]
which has a pair of complex conjugate poles at $\pm\mathrm{i}$ and is analytic in the strip $\mathcal{S}_{\rho}$ for $\rho\in(0,1)$. Figure \ref{fig:Super1} shows the first-order derivative error curve, i.e., $(f-h_n){'}(x)$, and the errors at the superconvergence points $\{\tau_j\}_{j=0}^{n+1}$. Figure \ref{fig:Super2} shows the second-order derivative error curve, i.e., $(f-h_n){''}(x)$, and the errors at the superconvergence points $\{\eta_j\}_{j=0}^{n+2}$. We see that the errors at the superconvergence points are significantly smaller than the maximum error.
\end{example}

\begin{figure}[htbp]
\centering
\includegraphics[width=0.75\textwidth,height=0.55\textwidth]{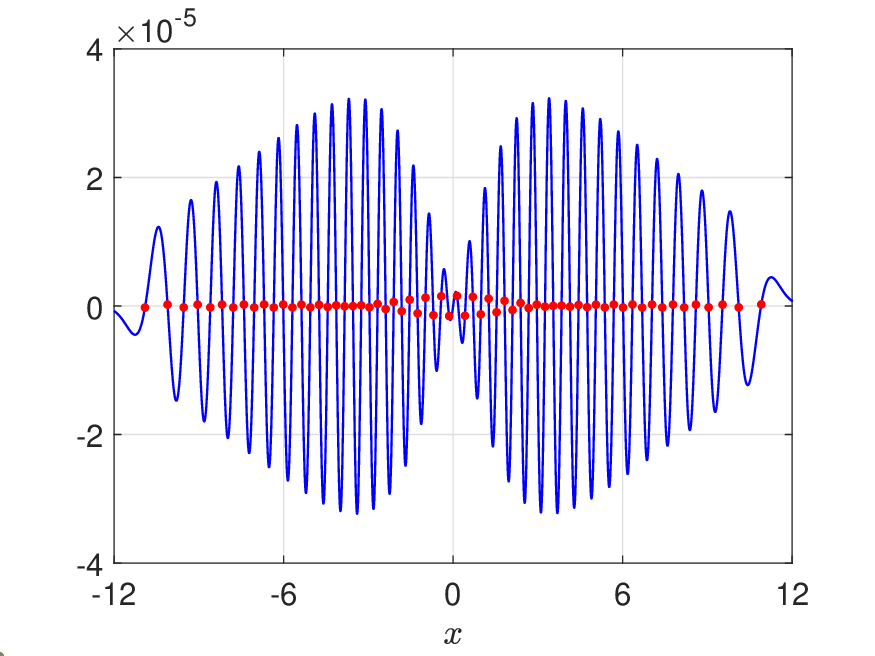}
\caption{Plot of $(f-h_n){'}(x)$ for $f(x)=\mathrm{e}^{-x^2}\cos(5x)$ and $n=62$ and the points are the errors at the superconvergence points $\{\tau_j\}_{j=0}^{n+1}$.}\label{fig:Super3}
\end{figure}

\begin{figure}[htbp]
\centering
\includegraphics[width=0.75\textwidth,height=0.55\textwidth]{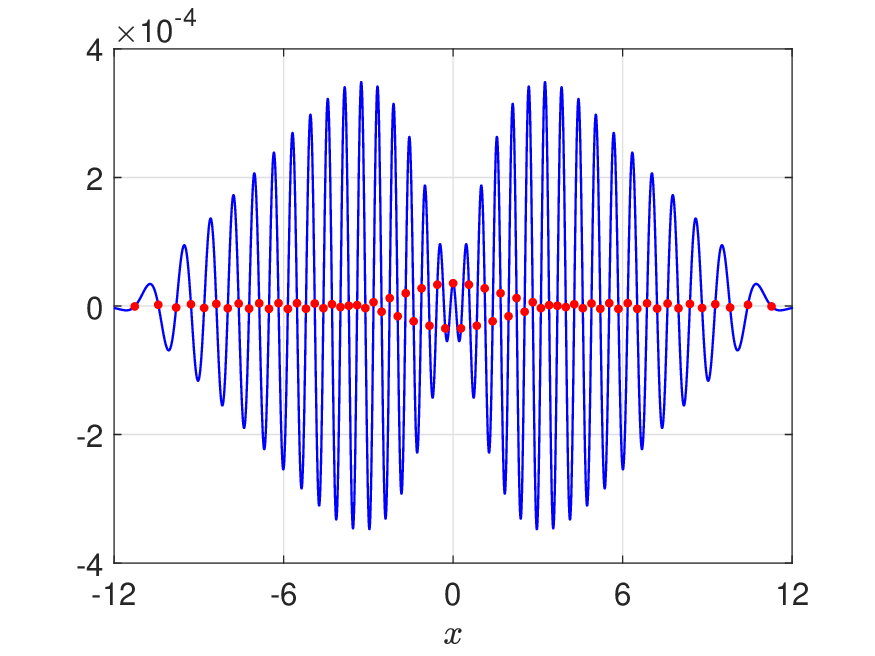}
\caption{Plot of $(f-h_n){''}(x)$ for $f(x)=\mathrm{e}^{-x^2}\cos(5x)$ and $n=62$ and the points are the errors at the superconvergence points $\{\eta_j\}_{j=0}^{n+2}$.}\label{fig:Super4}
\end{figure}

\begin{example}
We consider the wave packet function
\[
f(x) = \mathrm{e}^{-x^2}\cos(5x),
\]
which is analytic in the strip $\mathcal{S}_{\rho}$ for any $\rho\in(0,\infty)$. Figure \ref{fig:Super3} shows the first-order derivative error curve, i.e., $(f-h_n){'}(x)$, and the errors at the superconvergence points $\{\tau_j\}_{j=0}^{n+1}$. Figure \ref{fig:Super4} shows the second-order derivative error curve, i.e., $(f-h_n){''}(x)$, and the errors at the superconvergence points $\{\eta_j\}_{j=0}^{n+2}$. As expected, we see that the errors at the superconvergence points are much smaller than the maximum error.
\end{example}

Before closing this section, we show the sharpness of the factor $n^{1/2}$ gained at the superconvergence points for the first- and second-order derivatives of $h_n$. We define the following two ratios
\[
\mathcal{R}_1(n) = \frac{\max|(f - h_n){'}(\tau_j)|}{\|(f-h_n){'}\|_{\infty}}, \quad j=0,\ldots,n+1,
\]
and
\[
\mathcal{R}_2(n) = \frac{\max|(f - h_n){''}(\eta_j)|}{\|(f-h_n){''}\|_{\infty}}, \quad j=0,\ldots,n+2.
\]
In Figure \ref{fig:MaxError} we plot the maximum errors at the superconvergence points and the maximum errors over the real line and $\mathcal{R}_1(n)$ and $\mathcal{R}_2(n)$ multiplied by $n^{1/2}$, respectively, as functions of $n$ for $f(x)=\mathrm{e}^{-x^2/2}/(x^2+1)$. Clearly, we see that the maximum errors at the superconvergence points are remarkably smaller than the corresponding maximum errors over the real line. Moreover, we also see that $n^{1/2}\mathcal{R}_1(n)$ and $n^{1/2}\mathcal{R}_2(n)$ tend to some finite constants as $n$ increases and thus the factor $n^{1/2}$ gained at the superconvergence points for the first- and second-order derivatives is sharp.

\section{Superconvergence points of Hermite spectral collocation method}\label{sec:HermSpecColl}
In this section we consider the superconvergence points of Hermite spectral collocation method. For ease of exposition, we restrict our attention to two ODE model problems.

We first consider the following second order linear ODE model   
\begin{equation}
\begin{aligned}\label{def:Model}
& u''(x) + (\alpha - x^2)u(x) = f(x),  \\
&\lim_{x\rightarrow\pm\infty} u(x) = 0,
\end{aligned}
\end{equation}
where $\alpha\in\mathbb{R}$ and $\alpha\neq1,3,5,\ldots$ and $x\in\mathbb{R}$. The Hermite spectral collocation method is achieved by finding $u_n\in\mathbb{H}_{n}$ such that
\begin{equation}\label{eq:HermSpec}
u''_n(x_j) + (\alpha - x_j^2)u_n(x_j) = f(x_j), \quad j=0,\ldots,n,
\end{equation}
where $\{x_j\}_{j=0}^{n}$ are the zeros of $\psi_{n+1}(x)$. If we write $u_n$ in the Lagrange form
\begin{equation}
u_n(x) = \sum_{j=0}^{n} u(x_j) \sigma_j(x), \quad  \sigma_j(x) = \frac{\psi_{n+1}(x)}{\psi'_{n+1}(x_j)(x-x_j)},
\end{equation}
then \eqref{eq:HermSpec} can be rewritten as the following linear system
\begin{equation}
(D + S ) \hat{u}_n = \hat{f}_n,
\end{equation}
where $\hat{u}_n=(u(x_0),\ldots,u(x_n))^T$ and $\hat{f}_n=(f(x_0),\ldots,f(x_n))^T$ and
\begin{align}
D = \left( \begin{array}{ccc}
         \sigma''_0(x_0) & \cdots & \sigma''_n(x_0) \\
         \vdots & \ddots & \vdots \\
         \sigma''_0(x_n) & \cdots & \sigma''_n(x_n) \\
        \end{array}
        \right) , \quad
        S = \left(
                      \begin{array}{ccc}
                        \alpha - x_0^2 &  &  \\
                         & \ddots &  \\
                         &  & \alpha - x_n^2 \\
                      \end{array}
                    \right).
\end{align}
Below we consider the superconvergence points of Hermite spectral collocation method \eqref{eq:HermSpec}. From \cite{Babuska1996,Zhang2012} we know that the superconvergence property of finite element and polynomial spectral methods may be narrowed down to the capability of a polynomial space to approximate polynomials of one order higher. Here we generalize this criterion to the superconvergence analysis of Hermite spectral method and consider the capability of the space $\mathbb{H}_{n}$ to approximate functions in $\mathbb{H}_{n+1}$.
\begin{theorem}\label{thm:Supercon}
If $u\in\mathbb{H}_{n+1}$, then the Hermite spectral collocation method \eqref{eq:HermSpec} superconverges at $\{x_j\}_{j=0}^{n}$ and the first- and second-order derivatives superconverge at the zeros of $\psi'_{n+1}(x)$ and $\psi''_{n+1}(x)$, respectively.
\end{theorem}
\begin{proof}
Since $u\in\mathbb{H}_{n+1}$ and $u_n\in\mathbb{H}_n$, we write them in the forms $u(x)=a_0\psi_0(x)+\cdots+a_{n+1}\psi_{n+1}(x)$ and $u_n(x)=\hat{a}_0\psi_0(x)+\cdots+\hat{a}_{n}\psi_n(x)$, respectively. By \eqref{eq:harm} we see that
\begin{align}
u''_n(x) + (\alpha - x^2) u_n(x) &= \sum_{k=0}^{n} \hat{a}_k \left( \psi''_k(x) + (\alpha - x^2) \psi_k(x)  \right) = \sum_{k=0}^{n} \hat{a}_k (\alpha - 2k - 1) \psi_k(x). \nonumber
\end{align}
Note that the last sum belongs to $\mathbb{H}_n$ and, by \eqref{eq:HermSpec}, it interpolates $f(x)$ at the points $\{x_j\}_{j=0}^{n}$. By \cite[Equation~(7.90)]{Shen2011} we deduce that
\[
\hat{a}_k = \frac{1}{\alpha-2k-1} \sum_{j=0}^{n} \hat{w}_j \psi_k(x_j) f(x_j), \quad
k=0,\ldots,n,
\]
where $\hat{w}_j = ((n+1)(\psi_n(x_j))^2)^{-1}$. On the other hand, by \eqref{eq:harm} again we see that
\begin{align}
u{''}(x) + (\alpha - x^2) u(x) = \sum_{k=0}^{n+1} a_k (\alpha - 2k - 1) \psi_k(x) = f(x). \nonumber
\end{align}
Combining the above two results we find that
\begin{align}
\hat{a}_k &= \frac{1}{\alpha - 2k - 1} \sum_{j=0}^{n} \hat{w}_j \psi_k(x_j) \left( \sum_{\ell=0}^{n+1} a_{\ell} (\alpha - 2\ell - 1) \psi_{\ell}(x_j) \right) \nonumber \\
&= \frac{1}{\alpha - 2k - 1}  \sum_{\ell=0}^{n+1} a_{\ell} (\alpha - 2\ell - 1) \left( \sum_{j=0}^{n} \hat{w}_j \psi_k(x_j) \psi_{\ell}(x_j) \right) \nonumber \\
&= a_k, \nonumber
\end{align}
where we have used the discrete orthogonality satisfied by Hermite functions, i.e.,
\[
\sum_{j=0}^{n} \hat{w}_j \psi_k(x_j) \psi_{\ell}(x_j) = \delta_{k,\ell},
\]
in the last step. Hence, we deduce immediately that $u(x)-u_n(x)=a_{n+1}\psi_{n+1}(x)$ and thus $u(x)-u_n(x)$ superconverges at the points $\{x_j\}_{j=0}^{n}$ and its first- and second-order derivatives superconverge at the zeros of $\psi'_{n+1}(x)$ and $\psi''_{n+1}(x)$ respectively. This ends the proof.
\end{proof}

\begin{figure}[htbp]
\centering
\includegraphics[width=0.7\textwidth,height=0.55\textwidth]{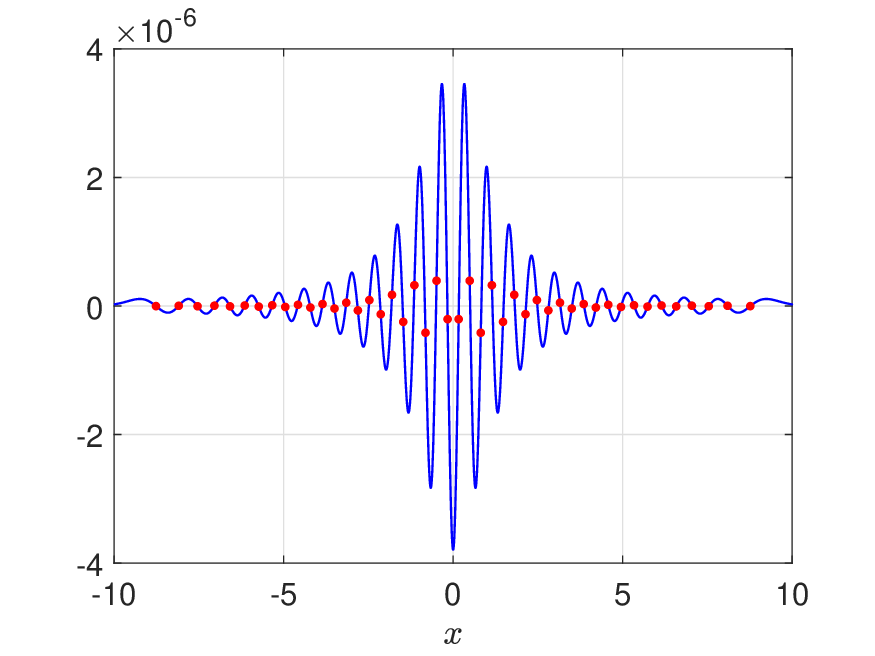}
\caption{Plot of $(u-u_n)(x)$ for $n=45$ and the points are the errors at the collocation points $\{x_j\}_{j=0}^{n}$. Here $\alpha=1/2$ and the exact solution is $u(x)=\mathrm{e}^{-x^2/2}/(x^2+2)$.}\label{fig:HermSpec1}
\end{figure}

\begin{figure}[htbp]
\centering
\includegraphics[width=0.7\textwidth,height=0.55\textwidth]{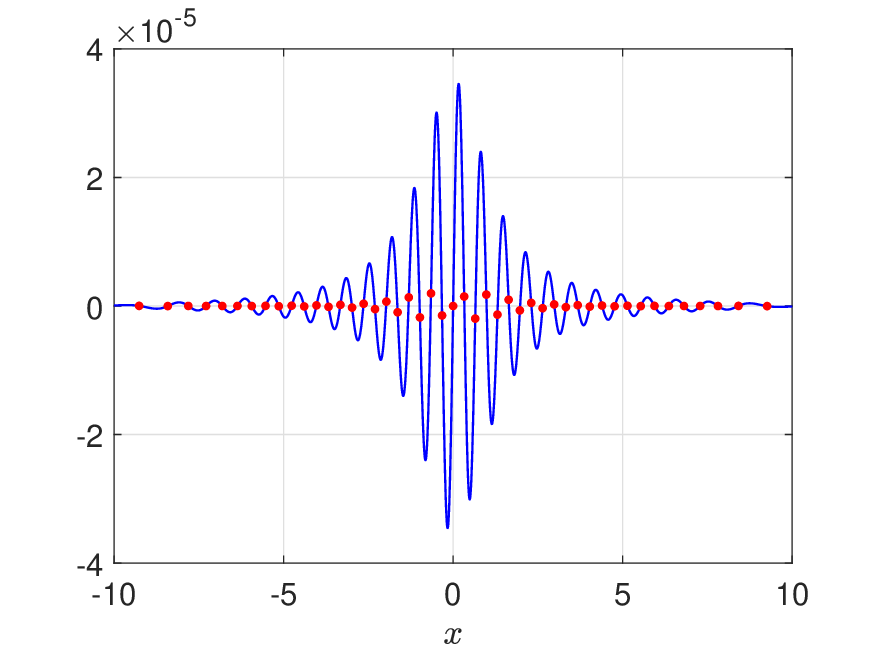}
\caption{Plot of $(u-u_n){'}(x)$ for $n=45$ and the points are the errors at the superconvergence points $\{\tau_j\}_{j=0}^{n+1}$. Here $\alpha=1/2$ and the exact solution is $u(x)=\mathrm{e}^{-x^2/2}/(x^2+2)$.}\label{fig:HermSpec2}
\end{figure}

To confirm our analysis, we consider the equation \eqref{def:Model} with $\alpha=1/2$ and $f(x)$ is chosen such that the exact solution is $u(x)=\mathrm{e}^{-x^2/2}/(x^2+2)$. In Figure \ref{fig:HermSpec1} we plot the pointwise error of the Hermite spectral collocation solution, i.e., $(u-u_n)(x)$, and the errors at the collocation points $\{x_j\}_{j=0}^{n}$. We see that the errors at the collocation points $\{x_j\}_{j=0}^{n}$ are significantly smaller than the maximum error. In Figure \ref{fig:HermSpec2} we plot the first-order derivative error of the Hermite spectral collocation solution, i.e., $(u-u_n){'}(x)$, and the errors at the superconvergence points $\{\tau_j\}_{j=0}^{n+1}$. We see that the errors at the superconvergence points $\{\tau_j\}_{j=0}^{n+1}$ are significantly smaller than the maximum error.

Next we consider the following second order linear ODE model \cite[Section~7.4.2]{Shen2011}
\begin{equation}
\begin{aligned}\label{def:Model2}
&{-}u{''}(x) + \alpha u(x) = f(x),  \\
&\lim_{x\rightarrow\pm\infty} u(x) = 0,
\end{aligned}
\end{equation}
The Hermite spectral collocation method for \eqref{def:Model2} is achieved by finding $u_n\in\mathbb{H}_{n}$ such that
\begin{equation}\label{eq:HermSpec2}
{-}u''_n(x_j) + \alpha u_n(x_j) = f(x_j), \quad j=0,\ldots,n,
\end{equation}
where $\{x_j\}_{j=0}^{n}$ are the zeros of $\psi_{n+1}(x)$. Equivalently, \eqref{eq:HermSpec2} can also be rewritten as the following linear system
\begin{equation}\label{eq:HermSpec3}
(-D + \alpha I ) \hat{u}_n = \hat{f}_n,
\end{equation}
where $D,\hat{u}_n, \hat{f}_n$ are defined as above and $I$ is the identity matrix.

\begin{remark}
From \cite[Theorem 4]{Weideman1992} we know that the eigenvalues of $D$ are real, negative, distinct and given by $\lambda=-\mu_j^2$, where $\mu_j$ are the $\lfloor (n+1)/2 \rfloor$ positive roots of $\psi_{n+1}(x)$, together with the interlacing $\lfloor (n+2)/2 \rfloor$ positive roots of $\psi'_{n+1}(x)$. Thus, to ensure that the Hermite spectral collocation method \eqref{eq:HermSpec3} admits a unique solution, we always assume that $\alpha\neq -\mu_j^2$ in the subsequent analysis.
\end{remark}
Below we show that the superconvergence results of the Hermite spectral collocation method for \eqref{def:Model2} are the same as that of \eqref{def:Model}.
\begin{theorem}\label{thm:Super2}
If $u\in\mathbb{H}_{n+1}$, then the Hermite spectral collocation method \eqref{eq:HermSpec3} superconverges at $\{x_j\}_{j=0}^{n}$ and the first- and second-order derivatives superconverge at the zeros of $\psi'_{n+1}(x)$ and $\psi''_{n+1}(x)$, respectively.
\end{theorem}
\begin{proof}
Following the idea of the proof of Theorem \ref{thm:Supercon}, we first write $u$ and $u_n$ in the forms
\[
u(x) = \sum_{k=0}^{n+1} a_k \psi_k(x), \quad  u_n(x) = \sum_{k=0}^{n} \hat{a}_k \psi_k(x).
\]
By \eqref{eq:harm} we see that
\begin{align}
{-}u''_n(x) + \alpha u_n(x) &= \sum_{k=0}^{n} \hat{a}_k \left({-}\psi''_k(x) + \alpha \psi_k(x)  \right) = \sum_{k=0}^{n} \hat{a}_k \left({-}x^2 + \alpha + 2k + 1\right) \psi_k(x), \nonumber
\end{align}
and by \eqref{eq:HermSpec2},
\begin{align}\label{eq:ColEqn1}
f(x_j) = \sum_{k=0}^{n} \hat{a}_k \left(-x_j^2 + \alpha + 2k + 1\right) \psi_k(x_j), \quad j=0,\ldots,n.
\end{align}
On the other hand, by \eqref{eq:harm} again we have
\begin{align}
-u{''}(x) + \alpha u(x) &= \sum_{k=0}^{n+1} a_k \left({-}\psi''_k(x) + \alpha \psi_k(x)  \right) = \sum_{k=0}^{n+1} a_k \left({-}x^2 + \alpha + 2k + 1\right) \psi_k(x), \nonumber
\end{align}
and thus
\begin{align}\label{eq:ColEqn2}
f(x_j) = \sum_{k=0}^{n+1} a_k \left({-}x_j^2 + \alpha + 2k + 1\right) \psi_k(x_j) = \sum_{k=0}^{n} a_k \left({-}x_j^2 + \alpha + 2k + 1\right) \psi_k(x_j).
\end{align}
Comparing \eqref{eq:ColEqn1} and \eqref{eq:ColEqn2} and noting \eqref{eq:HermSpec3} admits a unique solution, we see immediately that $a_k=\hat{a}_k$ for $k=0,\ldots,n$ and thus $u(x)-u_n(x)=a_{n+1}\psi_{n+1}(x)$. Hence, the desired results follow and this ends the proof.
\end{proof}

\begin{figure}[htbp]
\centering
\includegraphics[width=0.7\textwidth,height=0.55\textwidth]{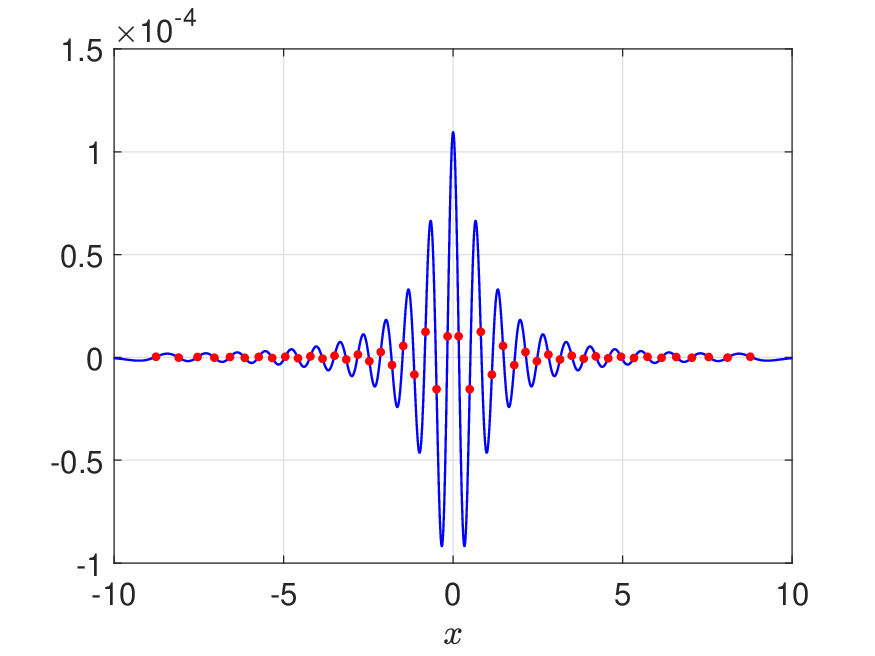}
\caption{Plot of $(u-u_n)(x)$ for $n=45$ and the points are the errors at the points $\{x_j\}_{j=0}^{n}$. Here $\alpha=2$ and the exact solution is $u(x)=\mathrm{e}^{-x^2}\ln(x^2+1)$.}\label{fig:HermSpec3}
\end{figure}

\begin{figure}[htbp]
\centering
\includegraphics[width=0.7\textwidth,height=0.55\textwidth]{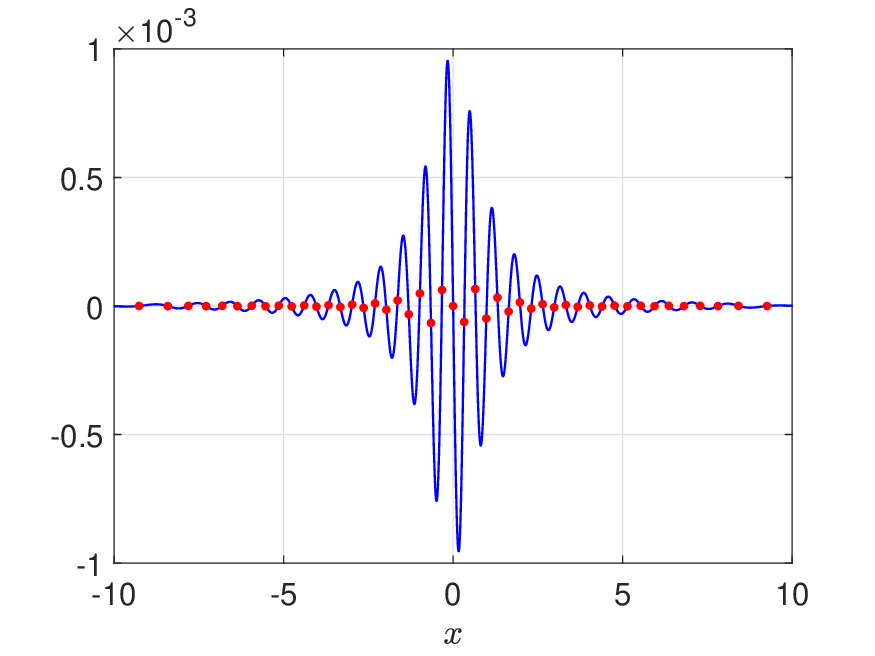}
\caption{Plot of $(u-u_n){'}(x)$ for $n=45$ and the points are the errors at the superconvergence points $\{\tau_j\}_{j=0}^{n+1}$. Here $\alpha=2$ and the exact solution is $u(x)=\mathrm{e}^{-x^2}\ln(x^2+1)$.}\label{fig:HermSpec4}
\end{figure}

We consider \eqref{def:Model2} with $\alpha=2$ and $f(x)$ is chosen such that the exact solution is $u(x)=\mathrm{e}^{-x^2}\ln(x^2+1)$. In Figure \ref{fig:HermSpec3} we plot the error curve of $(u-u_n)(x)$ and the errors at the collocation points $\{x_j\}_{j=0}^{n}$. We see that the errors at the collocation points $\{x_j\}_{j=0}^{n}$ are significantly smaller than the maximum error. In Figure \ref{fig:HermSpec4} we plot the first-order derivative error of the Hermite spectral collocation solution, i.e., $(u-u_n){'}(x)$, and the errors at the superconvergence points $\{\tau_j\}_{j=0}^{n+1}$. We see that the errors at the superconvergence points are significantly smaller than the maximum error.

\section{Post-processing of Hermite spectral collocation method}\label{sec:Post}
In this section we consider to extend the superconvergence analysis to post-processing of Hermite spectral collocation method.

\begin{figure}[htbp]
\centering
\includegraphics[width=0.7\textwidth,height=0.55\textwidth]{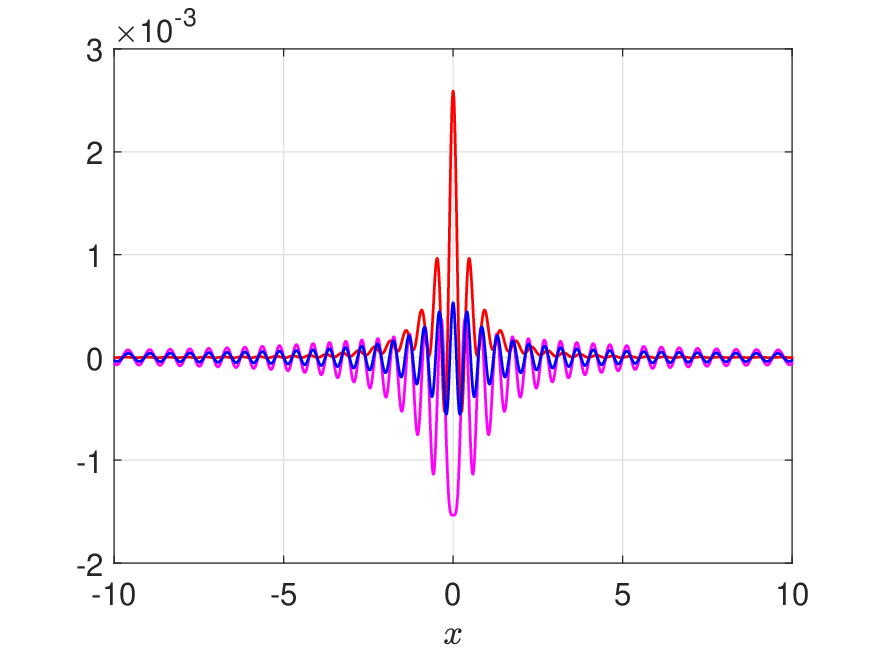}
\caption{Plot of $(u-u_n)(x)$ (magenta), $(u-u_{n+1})(x)$ (red) and $(u-\varphi_m)(x)$ (blue) for $n=90$ and $m=91$. Here $\alpha=1$ and the exact solution is $u(x)=(\mathrm{e}^{-(x-1)^2} + \mathrm{e}^{-(x+1)^2})/(4x^2+1)$.}\label{fig:HermPost1}
\end{figure}

\begin{figure}[htbp]
\centering
\includegraphics[width=0.7\textwidth,height=0.55\textwidth]{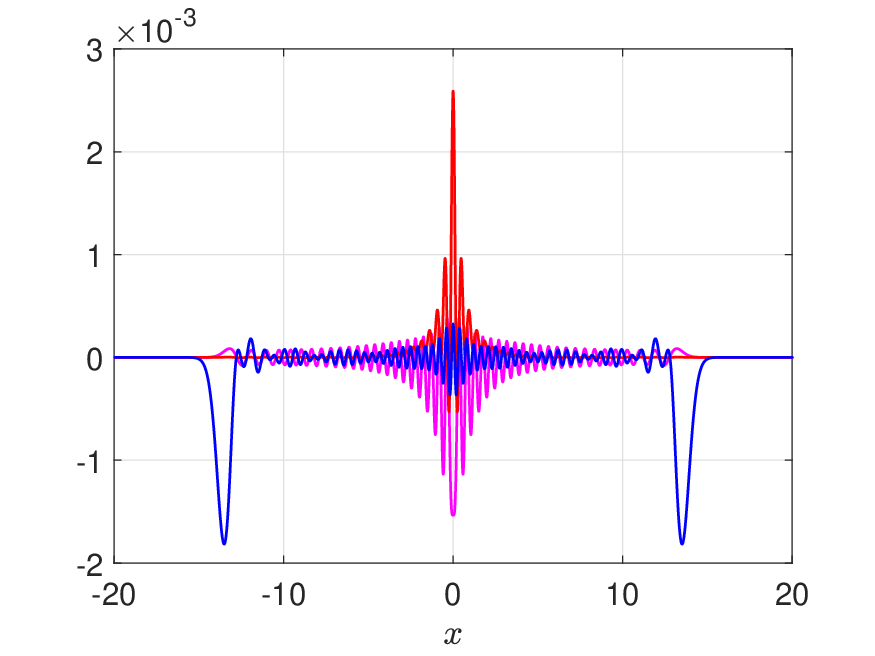}
\caption{Plot of $(u-u_n)(x)$ (magenta), $(u-u_{n+1})(x)$ (red) and $(u-\varphi_m)(x)$ (blue) for $n=90$ and $m=101$. Here $\alpha=1$ and the exact solution is $u(x)=(\mathrm{e}^{-(x-1)^2} + \mathrm{e}^{-(x+1)^2})/(4x^2+1)$.}\label{fig:HermPost2}
\end{figure}

\begin{figure}[htbp]
\centering
\includegraphics[width=0.7\textwidth,height=0.55\textwidth]{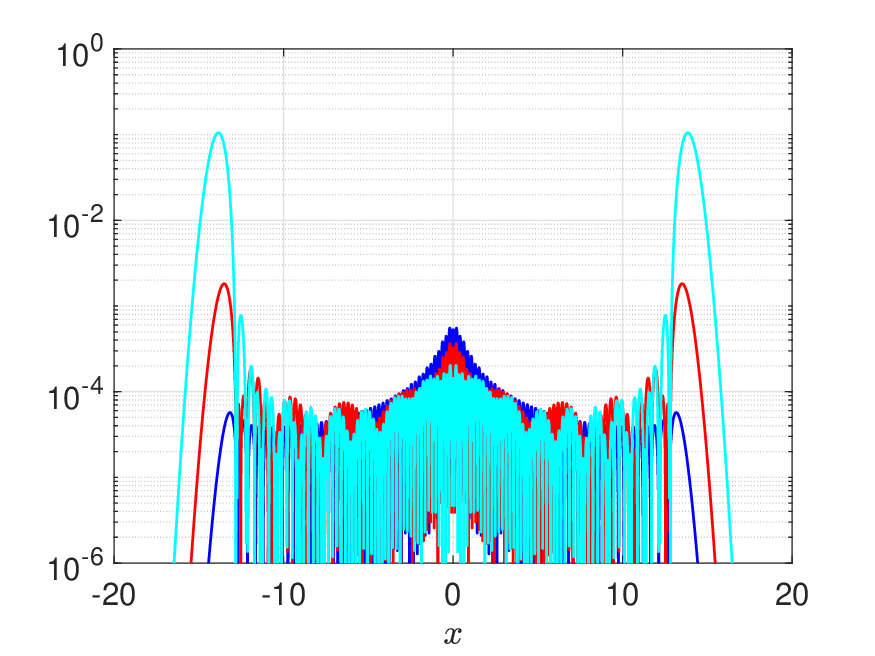}
\caption{Plot of $|(u-\varphi_m)(x)|$ for $m=91$ (blue), $m=101$ (red) and $m=111$ (cyan). Here $n=90$, $\alpha=1$ and the exact solution is $u(x)=(\mathrm{e}^{-(x-1)^2} + \mathrm{e}^{-(x+1)^2})/(4x^2+1)$.}\label{fig:HermPost3}
\end{figure}

We consider the ODE model \eqref{def:Model2} and denote by $u_n(x)$ the solution of Hermite spectral collocation method in \eqref{eq:HermSpec2}. Suppose that we have computed two solutions $u_{n}(x)$ and $u_{n+1}(x)$. From Theorem \ref{thm:Super2} we know that those two solutions superconverge at the zeros of $\psi_{n+1}(x)$ and $\psi_{n+2}(x)$, respectively. Now we consider to seek a new approximation in the space $\mathbb{H}_m$ ($m\le 2n+1$) by minimizing the errors at the superconvergence points of $u_{n}(x)$ and $u_{n+1}(x)$:
\begin{equation}
\varphi_m(x) = \mathrm{arg}\min_{\varphi\in\mathbb{H}_m} \left( \sum_{j=0}^{n} \left| \varphi(x_j) - u_{n}(x_j) \right|^2 + \sum_{j=0}^{n+1} \left| \varphi(y_j) - u_{n+1}(y_j) \right|^2 \right),
\end{equation}
where $\{x_j\}_{j=0}^{n}$ are the zeros of $\psi_{n+1}(x)$ and $\{y_j\}_{j=0}^{n+1}$ are the zeros of $\psi_{n+2}(x)$. If we set $\varphi_m(x)=a_0\psi_0(x)+\cdots+a_m\psi_m(x)$, then the above problem is equivalent to finding the least squares solution of the following problem:
\begin{align}\label{eq:LSProb}
\mathbf{a} = \min_{\mathbf{x}\in\mathbb{R}^{m+1}} \left\|\left(
  \begin{array}{c}
    \mathbf{A}_{1} \\
    \mathbf{A}_{2} \\
  \end{array}
\right) \mathbf{x} - \left(
                       \begin{array}{c}
                         \mathbf{u}_{1} \\
                         \mathbf{u}_{2} \\
                       \end{array}
                     \right) \right\|,
\end{align}
where $\|\cdot\|$ is the 2-norm for vectors and $\mathbf{a}=(a_0,\ldots,a_m)^T$ and
\begin{align}
\mathbf{A}_1 = \left(
           \begin{array}{ccc}
             \psi_0(x_0) & \cdots & \psi_m(x_0) \\
             \vdots & \ddots & \vdots \\
             \psi_0(x_n) & \cdots & \psi_m(x_n) \\
           \end{array}
         \right),
         \quad
\mathbf{A}_2 = \left(
           \begin{array}{ccc}
             \psi_0(y_0)   & \cdots & \psi_m(y_0) \\
             \vdots & \ddots & \vdots \\
             \psi_0(y_{n+1}) & \cdots & \psi_m(y_{n+1}) \\
           \end{array}
         \right), \nonumber
\end{align}
and
\begin{align}
\mathbf{u}_{1} = \left(
        \begin{array}{c}
          u_n(x_0) \\
          \vdots \\
          u_n(x_n) \\
        \end{array}
      \right), \quad
\mathbf{u}_{2} = \left(
        \begin{array}{c}
          u_{n+1}(y_{0}) \\
          \vdots \\
          u_{n+1}(y_{n+1}) \\
        \end{array}
      \right). \nonumber
\end{align}
In the following we show the performance of the post-processed solution $\varphi_m(x)$. We consider \eqref{def:Model2} with $\alpha=1$ and $f(x)$ is chosen such that the exact solution is $u(x)=(\mathrm{e}^{-(x-1)^2} + \mathrm{e}^{-(x+1)^2})/(4x^2+1)$. In Figure \ref{fig:HermPost1} we plot the errors of of two Hermite spectral collocation solutions $u_n(x)$ and $u_{n+1}(x)$ with $n=90$ and the post-processing solution $\varphi_m(x)$ with $m=91$. As expected, we see that the accuracy of the post-processed solution is better than that the two spectral collocation solutions. We further ask the question: Is it possible to improve the accuracy of the post-processed solution by increasing $m$? In Figure \ref{fig:HermPost2} we plot the errors of two Hermite spectral collocation solutions $u_n(x)$ and $u_{n+1}(x)$ with $n=90$ again, and the post-processed solution $\varphi_m(x)$ with $m=101$. We see that the accuracy of the post-processed solution $\varphi_m(x)$ can be improved further in some neighborhood of the origin, but deteriorates in neighborhoods of the largest and smallest zeros of $\psi_{n+2}(x)$ (i.e., $y_0$ and $y_{n+1}$). To show this more clearly, we further plot the absolute errors of $\varphi_m(x)$ for three values of $m$ in Figure \ref{fig:HermPost3} and we can see that the errors in neighborhoods of $y_0$ and $y_{n+1}$ grow very quickly when $m$ increases from $91$ to $111$. Indeed, the degradation in performance of $\varphi_m(x)$ with larger $m$ is due to overfitting, which is one of the most common problems in regression analysis and machine learning, and therefore we recommend choosing $m$ which is close to $n+1$ when computing $\varphi_m(x)$.

\section{Conclusion}\label{sec:Conclusion}
In this work, we have analyzed the superconvergence property of Hermite spectral interpolation. Based on the contour integral representation of Hermite spectral interpolation recently developed in \cite{Wang2025a}, we have identified the superconvergence points where the first- and second-order derivatives of Hermite spectral interpolation converge at a faster rate. We have further extended the superconvergence analysis to the Hermite spectral collocation method for two ODE models and derived similar superconvergence results. Moreover, we have also extended the superconvergence analysis to the post-processing of the Hermite spectral collocation method and showed that a more accurate approximation can be derived using the values of Hermite spectral collocation methods at the superconvergence points. Numerical results have been provided to confirm our analysis.

This work can be further extended along several directions, including Hermite spectral method for time-dependent PDEs and superconvergence points of Laguerre spectral interpolation. We leave these issues for future research.

\bigskip

\noindent {\bf Funding} The first author was partially supported by the National Natural Science Foundation of China under grant number 12371367 and the Hubei Provincial Natural Science Foundation of China under grant number 2023AFA083 and the fundamental research funds for the central universities under grant number 2025BRSXA003.

\bigskip

\noindent {\bf Data Availability} No data was used in this article.

\section*{Declarations}
{\bf Competing Interests} The authors declare that they have no conflict of interest.

\end{document}